\magnification=\magstep1          % |
\hsize=16truecm                   % | don't change these values !!
\vsize=23.5truecm                 % |
%

%  set \Finalout at the top of your paper, if you don't like
%  overfullrules!
\parindent=6mm
\normallineskiplimit=.99pt
\catcode`\@=11
%
%fonts (preloaded in PLAIN and AMS fonts)
%roman:
 \font\ninrm=cmr9 \font\egtrm=cmr8 \font\sixrm=cmr6
%math italic
 \font\nini=cmmi9 \skewchar\nini='177
 \font\egti=cmmi8 \skewchar\egti='177
 \font\sixi=cmmi6 \skewchar\sixi='177
%symbols
 \font\ninsy=cmsy9 \skewchar\ninsy='60
 \font\egtsy=cmsy8 \skewchar\egtsy='60
 \font\sixsy=cmsy6 \skewchar\sixsy='60
 \font\ninex=cmex9 \font\egtex=cmex8
%bold
 \font\ninbf=cmbx9 \font\egtbf=cmbx8 \font\sixbf=cmbx6
%text italic
 \font\ninit=cmti9 \font\egtit=cmti8
%slanted
 \font\ninsl=cmsl9 \font\egtsl=cmsl8
%caps and small caps
 \font\tensc=cmcsc10 \font\ninsc=cmcsc9 \font\egtsc=cmcsc8
%sans serif
 \font\tensf=cmss10 \font\ninsf=cmss9 \font\egtsf=cmss8
%euler and blackboard-bold
% \font\teneum=eufm10 \font\seveneum=eufm7 \font\fiveeum=eufm5   % EUL10
% \font\nineum=eufm9 \font\sixeum=eufm6                          % EUL9
% \font\tenmsb=msbm10 \font\sevenmsb=msbm7 \font\fivemsb=msbm5   % BBB10
% \font\ninmsb=msbm9 \font\sixmsb=msbm6                          % BBB9
%title
 \font\titbf=cmbx10 scaled \magstep2
 \font\titsyt=cmbsy10 scaled \magstep2 \skewchar\titsyt='60
 \font\titsys=cmbsy7 scaled \magstep2 \skewchar\titsys='60
 \font\titsyss=cmbsy5 scaled \magstep2 \skewchar\titsyss='60
 \font\titmit=cmmib10 scaled \magstep2 \skewchar\titmit='177
 \font\titmis=cmmib7 scaled \magstep2 \skewchar\titmis='177
 \font\titmiss=cmmib5 scaled \magstep2 \skewchar\titmiss='177
 \font\titit=cmbxti10 scaled \magstep2
% \font\titeut=eufb10 scaled \magstep2                           % EUL14
% \font\titeus=eufb7 scaled \magstep2                            % EUL14
% \font\titeuss=eufb5 scaled \magstep2                           % EUL14
% \font\titmsbt=msbm10 scaled \magstep2                          % BBB14
% \font\titmsbs=msbm7 scaled \magstep2                           % BBB14
% \font\titmsbss=msbm5 scaled \magstep2                          % BBB14
%sectionheadings
 \font\secbf=cmbx10 scaled \magstep1
 \font\secsyt=cmbsy10 scaled \magstep1 \skewchar\secsyt='60
 \font\secsys=cmbsy7 scaled \magstep1 \skewchar\secsys='60
 \font\secsyss=cmbsy5 scaled \magstep1 \skewchar\secsyss='60
 \font\secmit=cmmib10 scaled \magstep1 \skewchar\secmit='177
 \font\secmis=cmmib7 scaled \magstep1 \skewchar\secmis='177
 \font\secmiss=cmmib5 scaled \magstep1 \skewchar\secmiss='177
 \font\secit=cmbxti10 scaled \magstep1
% \font\seceut=eufb10 scaled \magstep1                           % EUL12
% \font\seceus=eufb7 scaled \magstep1                            % EUL12
% \font\seceuss=eufb5 scaled \magstep1                           % EUL12
% \font\secmsbt=msmb10 scaled \magstep1                          % BBB12
% \font\secmsbs=msbm7 scaled \magstep1                           % BBB12
% \font\secmsbss=msbm5 scaled \magstep1                          % BBB12
%subsectionheadings
 \font\ssecsyt=cmbsy10 \skewchar\ssecsyt='60
 \font\ssecsys=cmbsy7 \skewchar\ssecsys='60
 \font\ssecsyss=cmbsy5 \skewchar\ssecsyss='60
 \font\ssecmit=cmmib10 \skewchar\ssecmit='177
 \font\ssecmis=cmmib7 \skewchar\ssecmis='177
 \font\ssecmiss=cmmib5 \skewchar\ssecmiss='177
 \font\ssecit=cmbxti10
% \font\sseceut=eufb10 \font\sseceus=eufb7 \font\sseceuss=eufb5  % EUL10
%
%new families: euler and blackboard bold
%\newfam\eufam                                                   % EUL
%\def\frak{%                                                     % .
% \let\next\relax                                                % .
% \ifmmode
%  \let\next\frak@
% \else
%  \def\next{\errmessage{Use \string\frak\space in math mode only }}
% \fi
% \next}
%\def\frak@#1{{\frak@@{#1}}}                                     % .
%\def\frak@@#1{\fam\eufmfam#1}                                   % .
%\let\goth\frak                                                  % EUL
%
%\newfam\bbfam                                                   % BBB
%\def\bbb{%                                                      % .
% \let\next\relax                                                % .
% \ifmmode\let\next\bbb@
% \else\def\next{\errmessage{Use \string\bbb\space in math mode only}}
% \fi
% \next}                                                         % .
%\def\bbb@#1{{\bbb@@{#1}}}                                       % .
%\def\bbb@@#1{\fam\bbfam#1}                                      % BBB
%
%font sizes (adapted from MANMAC.TEX)
\newdimen\b@gsize
\def\b@g#1#2{{\hbox{$\left#2\vcenter to#1\b@gsize{}%
                    \right.\n@space$}}}
\def\big{\b@g\@ne}
\def\Big{\b@g{1.5}}
\def\bigg{\b@g\tw@}
\def\Bigg{\b@g{2.5}}
\def\tenpoint{%
 \textfont0=\tenrm \scriptfont0=\sevenrm \scriptscriptfont0=\fiverm%
 \def\rm{\fam0\tenrm}%
 \textfont1=\teni \scriptfont1=\seveni \scriptscriptfont1=\fivei%
 \textfont2=\tensy \scriptfont2=\sevensy \scriptscriptfont2=\fivesy%
 \textfont3=\tenex \scriptfont3=\tenex \scriptscriptfont3=\tenex%
 \def\it{\fam\itfam\tenit}\textfont\itfam=\tenit%
 \def\sl{\fam\slfam\tensl}\textfont\slfam=\tensl%
 \textfont\bffam=\tenbf \scriptfont\bffam=\sevenbf%
 \scriptscriptfont\bffam=\fivebf%
 \def\bf{\fam\bffam\tenbf}%
 \def\sc{\tensc}%
 \def\sf{\tensf}%
% \textfont\eufam=\teneum \scriptfont\eufam=\seveneum%           % EUL10
% \scriptscriptfont\eufam=\fiveeum%                              % EUL10
% \textfont\bbfam=\tenmsb \scriptfont\bbfam=\sevenmsb%           % BBB10
% \scriptscriptfont\bbfam=\fivemsb%                              % BBB10
 \normalbaselineskip=12pt%
 \setbox\strutbox=\hbox{\vrule height8.5pt depth3.5pt width\z@}%
 \abovedisplayskip12pt plus3pt minus9pt%
 \belowdisplayskip12pt plus3pt minus9pt%
 \abovedisplayshortskip\z@ plus3pt%
 \belowdisplayshortskip7pt plus3pt minus4pt%
 \normalbaselines\rm%
 \setbox\z@\vbox{\hbox{$($}\kern\z@}\b@gsize=1.2\ht\z@%
}
\def\ninepoint{%
 \textfont0=\ninrm \scriptfont0=\sixrm \scriptscriptfont0=\fiverm%
 \def\rm{\fam0\ninrm}%
 \textfont1=\nini \scriptfont1=\sixi \scriptscriptfont1=\fivei%
 \textfont2=\ninsy \scriptfont2=\sixsy \scriptscriptfont2=\fivesy%
 \textfont3=\ninex \scriptfont3=\ninex \scriptscriptfont3=\ninex%
 \def\it{\fam\itfam\ninit}\textfont\itfam=\ninit%
 \def\sl{\fam\slfam\ninsl}\textfont\slfam=\ninsl%
 \textfont\bffam=\ninbf \scriptfont\bffam=\sixbf%
 \scriptscriptfont\bffam=\fivebf%
 \def\bf{\fam\bffam\ninbf}%
 \def\sc{\ninsc}%
 \def\sf{\ninsf}%
% \textfont\eufam=\nineum \scriptfont\eufam=\sixeum%             % EUL9
% \scriptscriptfont\eufam=\fiveeum%                              % EUL9
% \textfont\bbfam=\ninmsb \scriptfont\bbfam=\sixmsb%             % BBB9
% \scriptscriptfont\bbfam=\fivemsb%                              % BBB9
 \normalbaselineskip=11pt%
 \setbox\strutbox=\hbox{\vrule height8pt depth3pt width\z@}%
 \abovedisplayskip11pt plus2.7pt minus8.1pt%
 \belowdisplayskip11pt plus2.7pt minus8.1pt%
 \abovedisplayshortskip\z@ plus2.7pt%
 \belowdisplayshortskip6.3pt plus2.7pt minus3.6pt%
 \normalbaselines\rm%
 \setbox\z@\vbox{\hbox{$($}\kern\z@}\b@gsize=1.2\ht\z@%
}
\def\@eightpoint{%
 \textfont0=\egtrm \scriptfont0=\sixrm \scriptscriptfont0=\fiverm%
 \def\rm{\fam0\egtrm}%
 \textfont1=\egti \scriptfont1=\sixi \scriptscriptfont1=\fivei%
 \textfont2=\egtsy \scriptfont2=\sixsy \scriptscriptfont2=\fivesy%
 \textfont3=\egtex \scriptfont3=\egtex \scriptscriptfont3=\egtex%
 \def\it{\fam\itfam\egtit}\textfont\itfam=\egtit%
 \def\sl{\fam\slfam\egtsl}\textfont\slfam=\egtsl%
 \textfont\bffam=\egtbf \scriptfont\bffam=\sixbf%
 \scriptscriptfont\bffam=\fivebf%
 \def\bf{\fam\bffam\egtbf}%
 \def\sc{\egtsc}%
 \def\sf{\egtsf}%
 \normalbaselineskip=9pt%
 \setbox\strutbox=\hbox{\vrule height7pt depth2pt width\z@}%
 \abovedisplayskip10pt plus2.4pt minus7.2pt%
 \belowdisplayskip10pt plus2.4pt minus7.2pt%
 \abovedisplayshortskip\z@ plus2.4pt%
 \belowdisplayshortskip5.6pt plus2.4pt minus3.2pt%
 \normalbaselines\rm%
 \setbox\z@\vbox{\hbox{$($}\kern\z@}\b@gsize=1.2\ht\z@%
}
\def\@titlefonts{%
 \textfont0=\titbf \def\rm{\titbf}%
 \textfont1=\titmit%
 \scriptfont1=\titmis \scriptscriptfont1=\titmiss%
 \textfont2=\titsyt%
 \scriptfont2=\titsys \scriptscriptfont2=\titsyss%
 \textfont\itfam=\titit \def\it{\titit}%
% \textfont\eufam=\titeut \scriptfont\eufam=\titeus%             % EUL14
% \scriptscriptfont\eufam=\titeuss%                              % EUL14
% \textfont\bbfam=\titmsbt \scriptfont\bbfam=\titmsbs%           % BBB14
% \scriptscriptfont\bbfam=\titmsbss%                             % BBB14
 \let\bf\rm \let\sf\rm \let\sc\rm \let\sl\it%
 \normalbaselineskip18pt%
 \setbox\strutbox\hbox{%
  \vrule height 12.6pt depth5.4pt width\z@}%
 \normalbaselines\rm%
 \setbox\z@\vbox{\hbox{$($}\kern\z@}\b@gsize=1.2\ht\z@%
}
\def\@sectionfonts{%
 \textfont0=\secbf \def\rm{\secbf}%
 \textfont1=\secmit%
 \scriptfont1=\secmis \scriptscriptfont1=\secmiss%
 \textfont2=\secsyt%
 \scriptfont2=\secsys \scriptscriptfont2=\secsyss%
 \textfont\itfam=\secit \def\it{\secit}%
% \textfont\eufam=\seceut \scriptfont\eufam=\seceus%             % EUL12
% \scriptscriptfont\eufam=\seceuss%                              % EUL12
% \textfont\bbfam=\secmsbt \scriptfont\bbfam=\secmsbs%           % BBB12
% \scriptscriptfont\bbfam=\secmsbss%                             % BBB12
 \let\bf\rm \let\sf\rm \let\sc\rm \let\sl\it%
 \normalbaselineskip14.5pt%
 \setbox\strutbox\hbox{%
  \vrule height 10.1pt depth4.4pt width\z@}%
 \normalbaselines\rm%
 \setbox\z@\vbox{\hbox{$($}\kern\z@}\b@gsize=1.2\ht\z@%
}
\def\@subsecfonts{\tenpoint%
 \textfont0=\tenbf \def\rm{\tenbf}%
 \textfont1=\ssecmit%
 \scriptfont1=\ssecmis \scriptscriptfont1=\ssecmiss%
 \textfont2=\ssecsyt%
 \scriptfont2=\ssecsys \scriptscriptfont2=\ssecsyss%
 \textfont\itfam=\ssecit \def\it{\ssecit}%
% \textfont\eufam=\sseceut \scriptfont\eufam=\sseceus%           % EUL10
% \scriptscriptfont\eufam=\sseceuss%                             % EUL10
 \let\bf\rm \let\sf\rm \let\sc\rm \let\sl\it%
 \normalbaselines\rm%
}
\def\vfootnote#1{%
  \insert\footins\bgroup%
  \interlinepenalty\interfootnotelinepenalty%
  \ninepoint%
  \splittopskip\ht\strutbox \splitmaxdepth\dp\strutbox%
  \floatingpenalty\@MM%
  \leftskip\z@skip \rightskip\z@skip \spaceskip\z@skip \xspaceskip\z@skip%
  \textindent{#1}\footstrut\futurelet\next\fo@t}
\def\makefootline{\baselineskip18pt\line{\the\footline}}
\newif\ifSec
\outer\def\section#1\par{\@beginsection{#1}}
\def\@beginsection#1{%
 \tenpoint%
 \vskip\z@ plus.1\vsize\penalty-250%
 \vskip\z@ plus-.1\vsize\vskip1.75\bigskipamount%
 \message{#1}\leftline{\noindent\@sectionfonts#1}%
 \nobreak\bigskip\noindent%
 \global\Sectrue
 \everypar{\ifSec\global\Secfalse\everypar{}\fi}}%
\outer\def\references{\@beginsection{References}\parindent20pt\ninepoint}%
\let\ref\item%
\def\subsection#1{\tenpoint
 \ifSec\else\bigskip\fi
 {\@subsecfonts#1}}%
\long\def\proclaim#1{%
 \medbreak
 \bf#1\enspace\it\ignorespaces}
\long\def\statement#1{%
 \medbreak
 \bf#1\enspace\rm\ignorespaces}
\def\endproclaim{%
 \ifdim\lastskip<\medskipamount
  \removelastskip\penalty55\medskip\fi
 \rm}
\let\endstatement\endproclaim
\def\qedrule{\hbox{\vrule height1.4ex depth0pt width1ex}}
\def\qed{\unskip\nobreak\quad\qedrule\medbreak}
\newskip\figskipamount \figskipamount\bigskipamount
\newbox\@caption
\def\figskip{\vskip\figskipamount}
\long\def\figure#1#2{%
 \setbox\@caption\hbox{\ninepoint\rm\ignorespaces#2}%
 \ifhmode%
  \vadjust{%
   \figskip%
   \line{\hfill\vbox to #1{\vfill}}%
   \ifvoid\@caption\else%
    \vskip\medskipamount%
    \ifdim \wd\@caption>\hsize%
     \noindent{\ninepoint\rm\ignorespaces#2\par}%
     \setbox\@caption\hbox{}%
    \else%
     \line{\hfil\box\@caption\hfil}%
    \fi%
   \fi%
   \figskip}%
 \else \ifvmode
  \figskip
  \line{\hfill\vbox to #1{\vfill}}
  \ifvoid\@caption\else
   \medskip
   \ifdim \wd\@caption>\hsize
    \noindent{\ninepoint\rm\ignorespaces#2\par}
    \setbox\@caption\hbox{}
   \else
    \line{\hfil\box\@caption\hfil}
   \fi
  \fi
  \figskip
 \fi\fi}
\newcount\@startpage
\@startpage=1
\def\@logo{\vtop{\@eightpoint
 \line{\strut\hfill}
 \line{\strut\hfill}
 \line{\strut\hfill}}}
\def\logo#1#2#3#4#5{\global\@startpage=#4 \global\pageno=#4
 \def\@logo{\vtop{\@eightpoint
 \line{\hfill Zeitschrift f\"ur Analysis und ihre Anwendungen}
 \line{\hfill Journal for Analysis and its Applications}
 \line{\hfill Volume #1 (#2), No.\ #3, #4--#5}}}}
\headline={%
 \ifnum\pageno=\@startpage
  \hfill
 \else
  \ifodd\pageno
   \strut\hfill{\ninepoint \sf \@runtitle}\qquad{\tenpoint\bf\folio}
  \else
   \strut{\tenpoint\bf\folio}\qquad{\ninepoint\sf\@runauthor}\hfill
  \fi
 \fi}
\footline={%
 \ifnum\pageno=\@startpage
  \hfill
 \else
  \hfill
 \fi}
\def\author#1{\def\@author{\ignorespaces#1}}
\def\@author{}
\def\runauthor#1{\def\@runauthor{\ignorespaces#1\ et.~al.}}
\def\@runauthor{\@author}
\def\title#1{\def\@title{\def\\{\hfill\egroup\line\bgroup\strut\hfill}#1}}
\def\@title{}
\def\runtitle#1{\def\@runtitle{\ignorespaces#1}}
\def\@runtitle{\@title}
\long\def\abstract#1\endabstract{\def\@abstract{\ignorespaces#1}}
\def\@abstract{}
\def\keywords#1{\def\@keywords{\ignorespaces#1}}
\def\@keywords{}
\def\primclass#1{%
 \def\@pclass{\ignorespaces#1}
 \global\let\classification\undefined
 \gdef\@class{}}
\def\secclasses#1{\def\@sclass{\ignorespaces#1}}
\def\classification#1{%
 \def\@class{\ignorespaces#1}
 \global\let\primclass\undefined
 \global\let\secclasses\undefined
 \gdef\@pclass{}\gdef\@sclass{}}
\newcount\addresscount
\addresscount=1
\newcount\addressnum
\def\address#1{%
 \expandafter\gdef\csname @address\number\addresscount \endcsname {#1}
 \global\advance\addresscount by 1}
\def\@addresses{
 \kern-10pt
 \addressnum=0
 \loop
  \ifnum\addressnum<\addresscount
   \advance\addressnum by 1
   \csname @address\number\addressnum \endcsname
   \ifnum\addressnum<\addresscount\hfill\break\fi
 \repeat
}
\long\def\maketitle{\begingroup
 \parindent0pt
\def\afootnote{%
  \insert\footins\bgroup%
  \interlinepenalty\interfootnotelinepenalty%
  \ninepoint%
  \splittopskip\ht\strutbox \splitmaxdepth\dp\strutbox%
  \floatingpenalty\@MM%
  \leftskip\z@skip \rightskip\z@skip \xspaceskip\z@skip%
  \footstrut\futurelet\next\fo@t}
 \vglue-12.5mm
 \@logo
 \vskip16.5pt
 \begingroup\@titlefonts
 \line\bgroup\strut\hfill\@title\hfill\egroup
 \endgroup
 \bigskip
 \centerline{\tenpoint\noindent\bf\hfill\@author\hfill}
 \vskip20mm
 \vbox{\ninepoint\noindent \bf Abstract.\ \rm\@abstract\par}
 \medskip
 \vbox{\ninepoint\noindent \bf Keywords:\ \it\@keywords\par}
 \smallskip
 \vbox{\ninepoint\noindent \bf AMS subject classification:\
  \ifx\classification\undefined
   \rm Primary\ \@pclass, secondary\ \@sclass\par
  \else
   \rm \@class
  \fi}
 \afootnote{\@addresses}
 \let\maketitle\relax
 \gdef\@title{} \gdef\@abstract{}
 \gdef\@keywords{}
 \gdef\@pclass{} \gdef\@sclass{} \gdef\@class{}
 \addressnum=0
 \loop
  \ifnum\addressnum<\addresscount
  \advance\addressnum by 1
  \expandafter\gdef\csname @address\number\addressnum \endcsname{}
 \repeat
 \endgroup
}
% switch to tenpoint
\catcode`@=12
\tenpoint

\author{M.~Frank and V.~M.~Manuilov}
\address{M.~Frank: University of Leipzig, Inst.~Math., D-04109 Leipzig,
F.~R.~G.}
\address{V.~M.~Manuilov: Moscow State University, Faculty Mech. Math., 117234
Moscow, Russia}
\title{Diagonalizing ''compact'' operators \\ on Hilbert W*-modules}

\runtitle{Diagonalizing ''compact'' operators}

\abstract
For W*-algebras $A$ and self-dual Hilbert $A$-modules $\cal M$ we show
that every self-adjoint, ''compact'' module operator on $\cal M$ is
diagonalizable. Some specific properties of the eigenvalues and of the
eigenvectors are described.
\endabstract


\primclass{47C15}
\secclasses{46L99, 46H25, 47A75.}

\maketitle

%\vfill \eject

%%% BEGINNING OF TEXT %%%
\smallskip \noindent
The goal of the present short note is to consider self-adjoint, ''compact''
module ope\-rators on self-dual Hilbert W*-modules (which can be supposed to
possess a countably generated W*-predual Hilbert W*-module, in general) with
respect to their diagonali\-zability. Some special properties of their
eigenvalues and eigenvectors are described.
A partial result in this direction was recently obtained by V.~M.~Manuilov
[10,11] who proved that every such operator on the standard countably generated
Hilbert
W*-module $l_2(A)$ over {\it finite} W*-algebras $A$ can be diagonalized on the
respective $A$-dual Hilbert $A$-module $l_2(A)'$. The same was shown to be
true for every self-adjoint bounded module operator on finitely generated
Hilbert C*-modules over general W*-algebras by R.~V.~Kadison
[5,6,7] and over commutative AW*-algebras by K.~Grove and G.~K.~Pe\-der\-sen
[4] sometimes earlier. M.~Frank has made an attempt to find a generalized
version of the Weyl-Berg theorem in the $l_2(A)'$ setting for some (abelian)
monotone complete C*-algebras which should satisfy an additional condition,
as well as a counterexample, cf.~[2].
Further results on generalizations of the Weyl-von Neumann-Berg theorem
can be found e.~g.~in papers of G.~J.~Murphy [12], S.~Zhang [15,16] and
H.~Lin [9] .

\noindent
We go on to investigate situations where non-finite W*-algebras appear as
coefficients of the special Hilbert W*-modules under consideration
(Proposition 5), and where arbitrary self-dual Hilbert W*-modules are
considered (Theorem 9). The applied techniques are rather different from that
in [10,11].
By the way, the results of V.~M.~Manuilov in [10,11] are obtained to be valid
for arbitrary self-adjoint, ''compact'' module operators on the self-dual
Hilbert $A$-module $l_2(A)'$ over finite W*-algebras (Pro\-position 3).
This generalizes [10] since in the situation of finite W*-algebras $A$ the
set of ''compact'' ope\-rators on $l_2(A)$ may be
definitely smaller than that on $l_2(A)'$, and the latter
may not contain all bounded module ope\-rators on $l_2(A)'$, in general.
We characterize the role of self-duality for getting adequate results in the
finite W*-case (Proposition 4). The final result of our investigations is
Theorem 9 describing the diagonalizability of ''compact'' operators on
self-dual Hilbert W*-modules in a great generality.

\medskip \noindent
We consider Hilbert W*-modules $\{ {\cal M}, \langle .,. \rangle \}$ over
general W*-algebras $A$ , i.~e.~(left) $A$-modules $\cal M$ together with
an $A$-valued inner product $\langle .,. \rangle: {\cal M} \times {\cal M}
\rightarrow A$ satisfying the conditions:
\medskip
(i) $\langle x,x \rangle \ge 0$ for every $x \in {\cal M}$.
\medskip
(ii) $\langle x,x \rangle =0$ if and only if $x=0$.
\medskip
(iii) $\langle x,y \rangle = \langle y,x \rangle^*$ for every $x,y \in {\cal
M}$.
\medskip
(iv) $\langle ax+by,z \rangle = a\langle x,z \rangle + b \langle y,z \rangle$
for
every $a,b \in A$, $x,y,z \in {\cal M}$.
\medskip
(v) ${\cal M}$ is complete with respect to the norm $\|x\|=\|\langle x,x
\rangle \|_A^{1/2}$.
\medskip \noindent
We always suppose, that the linear structures of the W*-algebra $A$ and of
the (left) $A$-module ${\cal M}$ are compatible, i.~e. $\lambda(ax)=(\lambda
a)x=
a(\lambda x)$ for every $\lambda \in {\bf C}$, $a \in A$, $x \in {\cal M}$.
Let us denote the $A$-dual Banach $A$-module of a Hilbert $A$-module
$\{ {\cal M}, \langle .,. \rangle \}$ by ${\cal M}'= \{ r:{\cal M}
\rightarrow A \, : \, r \, - \, {\rm {\it A}-linear} \; {\rm and} \; {\rm
bounded} \}$.

\noindent
Hilbert W*-modules have some very nice properties in contrast to general
Hilbert C*-modules: First of all, the $A$-valued inner product can always
be lifted to an $A$-valued inner product on the $A$-dual Hilbert $A$-module
${\cal M}'$ via the canonical embedding
of ${\cal M}$ into ${\cal M}'$, $x \rightarrow \langle .,x \rangle$, turning
${\cal M}'$ into a (left) self-dual Hilbert $A$-module, (${\cal M}'=
({\cal M}')'$). Moreover, one has the following criterion on self-duality:

\proclaim{Proposition 1.} {\rm [1, Thm.~3.2]}
\noindent
Let $A$ be a W*-algebra and
$\{ {\cal M}, \langle .,. \rangle \}$ be a Hilbert $A$-module. Then the
following conditions are equivalent:
\medskip
(i) ${\cal M}$ is self-dual.
\medskip
(ii) The unit ball of ${\cal M}$ is complete with respect to the topology
$\tau_1$
induced by the semi-norms $\{ f(\langle .,. \rangle^{1/2} \}$ on ${\cal M}$,
where $f$ runs over the normal states of $A$.
\medskip
(iii) The unit ball of ${\cal M}$ is complete with respect to the topology
$\tau_2$
induced by the linear functionals $\{ f(\langle .,x \rangle) \}$ on
${\cal M}$ where $f$ runs over the normal states of $A$ and $x$ runs over
${\cal M}$.
\endproclaim

\noindent
Furthermore, on self-dual Hilbert W*-modules every bounded module operator
has an adjoint, and the Banach algebra of all bounded module operators is
actually a W*-algebra. And last but not least, every bounded module operator
on a Hilbert W*-module $\{ {\cal M}, \langle .,. \rangle \}$ can be continued
to a unique bounded module operator on its $A$-dual Hilbert W*-module
${\cal M}'$ preserving the operator norm. (Cf.~[13].)

\noindent
We want to consider (self-adjoint,) ''compact'' module operators on Hilbert
W*-modules. By G.~G.~Kasparov [8] an $A$-linear bounded module operator
$K$ on a Hilbert $A$-module $\{ {\cal M}, \langle .,. \rangle \}$ is
''compact'' if it belongs to the norm-closed linear hull of the elementary
operators
\medskip
$\{ \theta_{x,y} \, : \, \theta_{x,y}(z) = \langle z,x \rangle y \, , \,
x,y \in {\cal M} \}$
\medskip \noindent
The set of all ''compact'' operators on ${\cal M}$ is denoted by $K_A({\cal
M})$.
By [13, Thm. 15.4.2] the C*-algebra $K_A({\cal M})$ is a two-sided ideal
of the set of all bounded, adjointable module operators $End_A^*({\cal M})$
on ${\cal M}$, and both these sets coincide if and only if ${\cal M}$ is
algebraically finitely generated as an $A$-module, (cf.~also [3,Appendix]).
This will be used below. Since we are
going to investigate single ''compact'' operators we make the useful
observation
that both the range of a given ''compact'' operator and the support of it are
Hilbert C*-modules generated by countably many elements with respect to
the norm topology or at least with respect to the $\tau_1$-topology,
(cf.~Proposition 1).
Hence, without loss of generality we can restrict our attention to countably
generated Hilbert W*-modules and their W*-dual Hilbert W*-modules.

\noindent
We are especially interested in the Hilbert W*-module
\medskip
$l_2(A) = \{ \{a_i:i \in {\bf N} \} \, : \, a_i \in A \, , \, \sum_i a_ia_i^*
\; {\rm converges} \; {\rm with} \; {\rm respect} \; {\rm to} \; \|.\|_A \}$
\smallskip
$\langle \{ a_i \},\{ b_i \} \rangle = \|.\|_A-\lim_{N \in {\bf N}}
\sum_{i=1}^N a_ib_i^* \, $ ,
\medskip  \noindent
and in its $A$-dual Hilbert W*-module
\medskip
$l_2(A)' = \left\{ \{a_i:i \in {\bf N} \} \, : \, a_i \in A \, , \, \sup_{N \in
{\bf N}}
\left\| \sum_{i=1}^N a_ia_i^* \right\| < \infty \right\}$
\smallskip
$\langle \{ a_i \},\{ b_i \} \rangle = {\rm w}^*-\lim_{N \in {\bf N}}
\sum_{i=1}^N a_ib_i^* \,$ ,
\medskip \noindent
because of G.~G.~Kasparov's stabilization theorem [8], stating that every
countably
generated Hilbert C*-module over a unital C*-algebra $A$ is a direct summand
of $l_2(A)$.

\statement{Definition 2.} Let $A$ be a W*-algebra and let $\{ {\cal M},
\langle .,. \rangle \}$ be a self-dual Hilbert $A$-module possessing a
countably generated Hilbert $A$-module as its $A$-predual. A bounded module
operator $T$ on ${\cal M}$ is {\it diagonalizable} if there exists a sequence
$\{x_i : i \in {\bf N} \}$ of non-trivial elements of ${\cal M}$ such that:
\smallskip \noindent
(i) $T(x_i) = \Lambda_ix_i$ for some elements $\Lambda_i \in A$, ($i \in {\bf
N}$),
\smallskip \noindent
(ii) The Hilbert $A$-submodule generated by the elements $\{x_i\}$ inside
${\cal M}$ has a trivial orthogonal complement.
\smallskip \noindent
(iii) The elements $\{ x_i :  i \in {\bf N} \}$ are pairwise orthogonal, and
the values $\{ p_i=\langle x_i,x_i \rangle : i \in {\bf N} \}$ are
projections in $A$.
\smallskip \noindent
(iv) The equality $\Lambda_ip_i=\Lambda_i$ holds for the projection $p_i$,
($i \in {\bf N}$).
\endstatement

\noindent
Note, that the eigenvalues and the eigenvectors are not uniquely determined
for the operator $T$ since $T(x)=\Lambda x$ implies $T(y)=\Lambda' y$ for
$\Lambda'=u\Lambda u^*$ and $y=ux$ for all unitaries $u \in A$. Moreover, the
eigenvalues of $T$ do not belong to the center of $A$, in general.
Consequently,
$T(ax)=a(\Lambda x) \not=\Lambda (ax)$, in general. That is, eigenvectors are
often not one-to-one related to $T$-invariant $A$-submodules of the Hilbert
$A$-module ${\cal M}$ under consideration.

\medskip \noindent
Now, we start our investigations decomposing $A$ into components of prescribed
type with respect to its direct integral representation. Denote by $p$ that
central projection of $A$ dividing $A$ into a finite part $pA$ and into
an infinite part $(1-p)A$. That means, that with respect to the direct integral
decomposition of $A$ the fibers are almost everywhere factors of type ${\rm
I}_n$,
$n < \infty$, or ${\rm II}_1$ inside $pA$ and almost everywhere factors of type
${\rm I}_\infty$ or ${\rm II}_\infty$ or ${\rm III}$ inside $(1-p)A$.
Analogously, the Hilbert $A$-module $l_2(A)$ decomposes into the direct sum
of two Hilbert $A$-modules $l_2(A)=l_2(pA) \oplus l_2((1-p)A)$, and every
bounded $A$-linear operator $T$ on $l_2(A)$ splits into the direct sum $T=
pT \oplus (1-p)T$, where each part acts only on the respective part of the
Hilbert $A$-module non-trivially and at the same time as an $A$-linear
operator.

%%%%%%%%%%%%%%%%%%%%%%%%%%%%%%%%%%%%%%%%%%%%%%%%%%%%%%%%%%%%%%%%%%%%%%%%%%%%%%

\noindent
Consequently, we can proceed considering W*-algebras $A$ of coefficients of
prescribed type. Our first goal is to revise the case of finite W*-algebras
investigated by V.~M.~Ma\-nui\-lov. There the
set $K_A(l_2(A)')$ does not coincide with the set $End_A(l_2(A)')$, and there
are always self-adjoint, bounded module operators $T$ on $l_2(A)'$ which can
not be diagonalized. For example, consider a self-adjoint, bounded linear
operator $T_o$ on a separable Hilbert space $H$ being non-diagonalizable,
(cf. Weyl's theorem). Using the decomposition $l_2(A)=\overline{A \otimes H}$
one obtains a self-adjoint, bounded module operator $T$ on $l_2(A)$ by the
formula $T(a \otimes h)=a \otimes T_o(h)$, ($a \in A$, $h \in H$). The
operator $T$ extends to an operator on $l_2(A)'$, and $T$ can not be
diagonalizable by assumption.
Surprisingly, V.~M.~Manuilov proved that every self-adjoint, ''compact''
operator on the standard countably generated Hilbert W*-module $l_2(A)$ over
finite W*-algebras $A$ can be diagonalized on the respective $A$-dual
Hilbert $A$-module $l_2(A)'$.
A careful study of his detailed proofs at [10], [11] brings to light that
for finite W*-algebras with infinite center the continuation of the
''compact'' operators to the respective $A$-dual Hilbert $A$-module is not
only a proof-technical necessity, but it is of principal character.
Self-duality has to be supposed to warrant the diagonalizability of all
self-adjoint ''compact'' module operators on ${\cal M} \subseteq l_2(A)'$ in
the finite case, and the key steps of the proof can be repeated one-to-one.
Consequently, we give the generalized formulation of V.~M.~Manuilov's
diagonalization theorem for the finite case, and we show additionally that
self-duality is an essential property of Hilbert W*-modules for finding a
(well-behaved) diagonalization of arbitrary ''compact'' module operators on
them, in general.

\proclaim{Proposition 3.} {\rm (cf.~[10], [11, Thm.4.1])} $\quad$
Let $A$ be a W*-algebra of finite type. Then eve\-ry self-adjoint, ''compact''
module operator $K$ on $l_2(A)'$ is diagonalizable. The sequence of
eigenvalues $\{ \Lambda_n : n \in {\bf N} \}$ of $K$ has the property
$\lim_{n \rightarrow \infty} \|\Lambda_n\|=0$.
The eigenvalues $\{ \Lambda_n : n \in {\bf N} \}$ of $K$ can be
chosen in such a way that $\Lambda_2 \leq \Lambda_4 \leq ... \leq 0 \leq ...
\leq \Lambda_3 \leq \Lambda_1$.
Moreover, for positive operators $K$ without kernel the eigenvectors
$\{ x_n : n \in {\bf N} \}$ may possess the property $\langle x_n,x_n \rangle
=1_A$, $(n \in {\bf N})$, in addition.
\endproclaim

\noindent
For the detailed (but extended) proof of this proposition see [11], (also
[10]). The proving technique relies mainly on spectral decomposition theory
of operators and on the center-valued trace on the finite W*-algebra $A$.

\proclaim{Proposition 4.} Let $A$ be a finite W*-algebra with infinite center.
Consider a Hilbert $A$-module ${\cal M}$ such that $l_{2}(A) \subset
{\cal M} \subseteq l_{2}(A)'$. Then the following two statements are
equivalent:
\smallskip \noindent
(i)  $\, {\cal M} = l_{2}(A)'$, i.e., ${\cal M}$ is self-dual.
\smallskip \noindent
(ii) Every positive ''compact'' module operator is diagonalizable inside
     ${\cal M}$ with comparable inside the positive cone of $A$ eigenvalues.
\endproclaim

\statement{Proof.} Note, that $l_2(A) \not= l_2(A)'$ by assumption.
Denote the standard orthonormal basis of $l_2(A)$ by
$\{ e_n : n \in {\bf N} \}$. If the center of $A$ is supposed to be infinite
dimensional then one finds a sequence of pairwise orthogonal non-trivial
projections $\{ p_n : n \in {\bf N} \} \in {\rm Z}(A)$ summing up to $1_A$ in
the sense of w*-convergence.
Fix a sequence of positive non-zero numbers $\{ \alpha_n : n \in {\bf N} \}$
monotonically converging to zero. The bounded module operator $K$ defined by
\smallskip
$K(e_1)=\left( \sum_{n=1}^\infty \alpha_n p_n e_n \right) \,$,
$\, K(e_j)=\alpha_j p_j e_1$ for $j \not= 1$
\smallskip \noindent
is a ''compact'' operator on $l_2(A)$. It easily continues to a ''compact''
operator on ${\cal M}$. As an exercise one checks that the eigenvalues of
$K$ are $\{ \alpha_1p_1, \alpha_2p_2, \cdots ,0, \cdots ,-\alpha_2p_2 \}$
(ordering by sign and norm and taking into account (iii) and (iv) of Definition
2),
and that the appropriate eigenvectors are
\smallskip
  $\{ p_1e_1, 1/ \sqrt{2} p_2(e_1+e_2), 1/ \sqrt{2} p_3(e_1+e_3), \cdots ,
  \{ (1_A-p_n)e_n : n \in {\bf N} \}, \cdots$
\smallskip
   $\cdots , 1/ \sqrt{2}  p_3(e_1-e_3), 1/ \sqrt{2} p_2(e_1-e_2)\}$.
\smallskip \noindent
The only way of making the eigenvalues comparable inside the positive cone
of $A$ preser\-ving Definition 2,(iii)-(iv) is to sum up the positive and the
negative eigenvalues separately. But, then the resulting eigenvector
\smallskip
$x=(1_A+(1+1/ \sqrt{2})(1_A-p_1), 1/ \sqrt{2} p_2, 1/ \sqrt{2}  p_3, \cdots,
1/ \sqrt{2} p_n, \cdots )$,
\smallskip \noindent
corresponding to the only positive eigenvalue $\sum_{n=1}^\infty \alpha_n p_n$
of $K$ does not belong to ${\cal M}$ any longer by assumption. This shows one
implication. The converse implication follows from Proposition 6.  \qed
\endstatement

%%%%%%%%%%%%%%%%%%%%%%%%%%%%%%%%%%%%%%%%%%%%%%%%%%%%%%%%%%%%%%%%%%%%%%%%%%%%%%

\noindent
The second big step is to investigate the case of infinite W*-algebras as
coefficients of the Hilbert W*-modules under consideration.
The result is characteristic for the situation in self-dual
Hilbert W*-modules over infinite W*-algebras, and quite different from that
in the finite W*-case, and elsemore, from the classical Hilbert space
situation.

\proclaim{Proposition 5.}
Let $A$ be a W*-algebra which possesses infinitely many pairwise orthogonal,
non-trivial projections $\{p_i : i \in {\bf N} \}$ equivalent to $1_A$
and summing up to $1_A$ in the sense of w*-convergence of the sum
$\sum_i p_i=1_A$. Then the Hilbert $A$-module $l_2(A)'$ equipped with its
standard $A$-valued inner product is isomorphic to the Hilbert $A$-module
$\{ A,\langle .,. \rangle_A \}$, where $\langle a,b \rangle_A=ab^*$.
\endproclaim

\statement{Proof.}
Suppose, the equivalence of the projections $\{p_i : i \in {\bf N} \}$ with
$1_A$ is realized by partial isometries $\{ u_i : i \in {\bf N} \} \in A$,
$p_i=u_iu_i^*$, $1_A=u_i^*u_i$. Then the mapping
\medskip
$S: l_2(A)' \rightarrow A \quad , \quad \{ a_i \} \rightarrow
{\rm w}^*-\lim ({\rm finite} \, \sum_i a_iu_i^*)$
\medskip  \noindent
with the inverse mapping
\medskip
$S^{-1}: A \rightarrow l_2(A)' \quad , \quad a \rightarrow \{ au_i \}$
\medskip \noindent
realizes the isomorphism of $l_2(A)'$ and of $A$ as Hilbert $A$-modules because
of Proposi\-tion~1. \qed
\endstatement

\proclaim{Corollary 6.}
Let $A$ be a W*-algebra of infinite type. Then every bounded
module operator $T$ on $l_2(A)'$ is diagonalizable, and the formula
\smallskip
$T(\{ a_i \}) = \langle \{ a_i \}, \{ u_i \} \rangle \Lambda_T \{ u_i \} $
\smallskip \noindent
holds for every $\{ a_i \} \in l_2(A)'$, some $\Lambda_T \in A$ and the
partial isometries $\{ u_i \} \in A$ described in the previous proof.
\endproclaim

\statement{Proof.}
Every W*-algebra of type ${\rm I}_\infty$, ${\rm II}_\infty$ or ${\rm III}$
possesses a set of partial isometries with properties described at
Proposition 3. The same is true for W*-algebras consisting only of parts of
these types. Now, translate the operator $T$ on $l_2(A)'$ to an operator
$STS^{-1}$ on $A$ and vice versa using Proposition 3, and take into account
that every bounded module operator on $A$ is a multiplication operator with
a concrete element (from the right). \qed
\endstatement

\proclaim{Corollary 7.}
Let $A$ be a W*-algebra without any fibers of type ${\rm I}_n$, $n < \infty$,
and ${\rm II}_1$ in its direct integral decomposition. Let ${\cal M}$ be a
self-dual Hilbert $A$-module possessing a countably generated $A$-predual
Hilbert $A$-module. Then every bounded
module operator $T$ on ${\cal M}$ is diagonalizable, and the formula
\smallskip
$T(x) = \langle x, u \rangle \Lambda_T u $
\smallskip \noindent
holds for every $x \in {\cal M}$, some $\Lambda_T \in A$ and a universal for
all $T$ eigenvector $u \in {\cal M}$.
\endproclaim

\statement{Proof.}
Since ${\cal M}$ has a countably Hilbert $A$-module as its
$A$-predual, ${\cal M}$ is a direct summand of the Hilbert $A$-module
$l_2(A)'$ by G.~G.~Kasparov's stabilization theorem ([8]). Hence, one has to
show the assertion for the self-dual Hilbert $A$-module $l_2(A)'$ only. For
further use denote the projection from $l_2(A)'$ onto ${\cal M}$ by $P$.
Consider the direct integral decomposition of $A$ over its center. Therein
every fiber is a W*-factor of type ${\rm I}_\infty$, ${\rm II}_\infty$
or ${\rm III}$ by assumption. Putting it into the $l_2(A)'$-context one
obtains that $A$ is isomorphic to $l_2(A)'$ either applying Corollary 6
fiberwise or constructing a suitable set of partial isometries $\{ u_i \}
\in A$ to make use of Proposition 5. Then in the same way as there the
diagonalization result turns out for arbitrary bounded module operators $T$
on $l_2(A)'$. To get the formula of Corollary 7 one has only to set
$u = P(\{ u_i \})$.   \qed
\endstatement

\statement{Remark.} Let $A$ be a ${\rm I}_\infty$-factor, for example.
Then there are self-adjoint elements $\Lambda_T$ in $A$ which can not be
diagonalized in a stronger sense. More precisely, there is no way of
representing any such operator as a sum $\sum \lambda_i P_i$ with $\lambda_i
\in {\bf C}={\rm Z}(A)$ and $P_i=P_i^*=P_i^2 \in A$ because of Weyl's theorem.
Therefore, the Corollaries 6 and 7 are the strongest results one could expect.
\endstatement

%%%%%%%%%%%%%%%%%%%%%%%%%%%%%%%%%%%%%%%%%%%%%%%%%%%%%%%%%%%%%%%%%%%%%%%%%%%%%%

\statement{Example 8.}
Consider the C*-algebra $A$ of all $2 \times 2$-matrices on the complex
numbers.
Set ${\cal M}=A^2$ with the usual $A$-valued inner product. Consider the
(''compact'') bounded module operator $K=\theta_{x,x}+\theta_{y,y}$ for
\medskip
$x= \left( \left( \matrix{1 & 0 \cr 0 & 3 \cr} \right) ,
          \left( \matrix{0 & 0 \cr 0 & 0 \cr} \right) \right)$,
$y= \left( \left( \matrix{0 & 0 \cr 0 & 0 \cr} \right) ,
           \left( \matrix{2 & 0 \cr 0 & 2 \cr} \right) \right)$.
\medskip \noindent
Eigenvectors of $K$ are $x,y \in A^2$, for example, and the respective
eigenvalues are
\medskip
$\Lambda_x=\left( \matrix{1 & 0 \cr 0 & 9 \cr} \right)$,
$\Lambda_y=\left( \matrix{4 & 0 \cr 0 & 4 \cr} \right)$.
\medskip \noindent
Remark, that one can not compare these eigenvalues as elements of the
positive cone of $A$. But, making another choice one arrives at that
situation described at Proposition 6:
\medskip
$x_1=  \left( \left( \matrix{1 & 0 \cr 0 & 0 \cr} \right) ,
          \left( \matrix{0 & 0 \cr 0 & 1 \cr} \right) \right)$,
$x_2=  \left( \left( \matrix{0 & 0 \cr 0 & 1 \cr} \right) ,
           \left( \matrix{1 & 0 \cr 0 & 0 \cr} \right) \right)$.
\medskip \noindent
Then the respective eigenvalues are
\medskip
$\Lambda_1=\left( \matrix{1 & 0 \cr 0 & 4 \cr} \right)$,
$\Lambda_2=\left( \matrix{4 & 0 \cr 0 & 9 \cr} \right)$.
\medskip \noindent
and they can be ordered as well as the eigenvectors $x_1,x_2$ are units.
\noindent
Last but not least, dropping out condition (iv) of Definition 2 one can
correlate
$K$-invariant submodules of ${\cal M}$ and eigenvectors of $K$. Simply, set
\medskip
$x_1= \left( \left( \matrix{1 & 0 \cr 1 & 0 \cr} \right) ,
          \left( \matrix{0 & 1 \cr 0 & 1 \cr} \right) \right)$,
$x_2= \left( \left( \matrix{0 & 1 \cr 0 & 1 \cr} \right) ,
           \left( \matrix{1 & 0 \cr 1 & 0 \cr} \right) \right)$.
\medskip \noindent
In this case the corresponding eigenvalues are
\medskip
$\Lambda_1=\left( \matrix{1 & 0 \cr 0 & 4 \cr} \right)$,
$\Lambda_2=\left( \matrix{4 & 0 \cr 0 & 9 \cr} \right)$.
\medskip \noindent
They can be ordered in the positive cone of $A$. But, the eigenvectors
corresponding to the $K$-invariant submodules of ${\cal M}$ can not be
selected to be units any longer.
\endstatement

\proclaim{Theorem 9.}
Let $A$ be a W*-algebra and ${\cal M}$ be a self-dual Hilbert $A$-module.
Then eve\-ry self-adjoint, ''compact'' module operator on ${\cal M}$
is diagonalizable. The sequence of eigenvalues $\{ \Lambda_n : n \in {\bf N}
\}$ of $K$ has the property $\lim_{n \rightarrow \infty} \|\Lambda_n\|=0$.
The eigenvalues $\{ \Lambda_n : n \in {\bf N} \}$ of $K$ can be
chosen in such a way that $\Lambda_2 \leq \Lambda_4 \leq ... \leq 0 \leq ...
\leq \Lambda_3 \leq \Lambda_1$, and that $\{\Lambda_n : n \geq 3 \}$ are
contained in the finite part of {\bf A}.
\endproclaim

\statement{Proof.} Both the $\tau_1$-closure of the range and of the support
of $K$ are self-dual Hilbert C*-modules possessing countably generated
$A$-predual Hilbert $A$-modules because of the ''compact''ness of $K$.
Hence, without loss of generality one can restrict the attention to self-dual
Hilbert W*-modules with countably generated W*-predual Hilbert
W*-modules formed as the $\tau_1$-completed direct sum of range and support
of $K$. As usual, on the kernel of $K$ one has the eigenvalue zero and a
suitable system of eigenvectors.

\noindent
Now, gluing Corollary 4 and Proposition 6 together the theorem turns out to
be true in the special case ${\cal M} = l_2(A)'$, (cf.~the remarks in the
beginning of the present note). The only loss may be that the eigenvectors
are not units, in general. Because of G.~G.~Kasparov's stabilization
theorem ([8]) ${\cal M}$ possesses an embedding into $l_2(A)'$
as a direct summand by assumption. Therefore, every self-adjoint, ''compact''
module operator $K$ on ${\cal M}$ can be continued to a unique such operator
on $l_2(A)'$ preserving the norm, simply applying the rule $K|_{{\cal M}^\bot}
=0$. The eigenvectors of this extension are elements of ${\cal M}$. The
Hilbert $A$-module ${\cal M}^\bot$ belongs to its kernel. This shows the
theorem. \qed
\endstatement

\statement{Remark.} For commutative AW*-algebras $A$ the statement of Theorem 9
is still true by [4]. The general AW*-case is open at present because of two
crucial unsolved problems in the AW*-theory:
\noindent
(i) Are the self-adjoint elements of ${\rm M}_n(A)$, $n \geq 2$, diagonalizable
for arbitrary (monotone complete) AW*-algebras $A$, or not?
\noindent
(ii) Does every finite (monotone complete) AW*-algebra possess a center-valued
trace, or not?
\endstatement

\statement{Remark.} One can extend the statement of Theorem 9 to the case of
normal, ''compact'' module operators dropping out only the ordering of the
eigenvalues. To see this note that for normal elements $K$ of the C*-algebra
$K_A({\cal M})$ there always exists a self-adjoint element $K' \in
K_A({\cal M})$ such that $K$ is contained in that C*-subalgebra of
$End_A({\cal M})$ generated by $K'$ and by the identity operator. Applying
functional calculus inside the W*-algebra $End_A({\cal M})$ the result turns
out. Beside this, it would be interesting to investigate some more general
variants
of the Weyl-von Neumann-Berg theorem for appropriate bounded module operators
on (self-dual) Hilbert W*-submodules over (finite) W*-algebras $A$ as those
obtained by  H.~Lin, G.~J.~Murphy and S.~Zhang.
\endstatement

\statement{Acknowledgement.} The second author thanks for the partial support
by the Russian Foundation for Fundamental Research (grant no. 94-01-00108a)
and by the International Science Foundation (grant no. MGM 000). The research
work was carried out during a stay at Leipzig which was part of a university
cooperation project financed by Deutscher Akademischer Austauschdienst.
We are very appreciated to the referees for their remarks on the first version
of
the present work.
\endstatement

\references

\item{[1]} Frank,~M.:
{\it Self-duality and C*-reflexivity of Hilbert C*-modules.}
Zeitschr.~Anal.~Anwendungen 9(1990), 165-176.
\smallskip

\item{[2]} Frank,~M.:
{\it Hilbert C*-modules over monotone complete C*-algebras and a Weyl-Berg
type theorem.} preprint 3/91, Universit\"at Leipzig, NTZ, 1991. To appear in
Math.~Nachrichten.
\smallskip

\item{[3]} Frank,~M.:
{\it Geometrical aspects of Hilbert C*-modules.} preprint 22/93, K{\o}benhavns
Universitet, Matematisk Institut, 1993.
\smallskip

\item{[4]} Grove,~K. and G.~K.~Pedersen:
{\it Diagonalizing matrices over C(X).} J.~Functional Analysis 59(1984), 64-89.
\smallskip

\item{[5]} Kadison,~R.~V.:
{\it Diagonalizing matrices over operator algebras.} Bull.~Amer.~Math.~Soc.
8(1983), 84-86.
\smallskip

\item{[6]} Kadison,~R.~V.:
{\it Diagonalizing matrices.} Amer.~J.~Math. 106(1984), 1451-1468.
\smallskip

\item{[7]} Kadison,~R.~V.:
{\it The Weyl theorem and block decompositions.}
In: Operator Algebras and Applications, v.~1, Cambridge: Cambridge University
Press 1988, pp.~109-117.
\smallskip

\item{[8]} Kasparov,~G.~G.:
{\it Hilbert $C^*$-modules: The theorems of Stinespring and Voiculescu.}
J.~Operator Theory 4(1980), 133-150.
\smallskip

\item{(9]} Lin,~H.:
{\it The generalized Weyl - von Neumann theorem and C*-algebra extensions.}
In: Algebraic Methods in Operator Theory. (eds.: R.~Curto and
P.~E.~T.~J{\o}rgensen), Birkh\"auser, Boston - Basel - Berlin, 1994.
\smallskip

\item{[10]} Manuilov,~V.~M.:
{\it Diagonalization of compact operators on Hilbert modules over W*-algebras
of finite type (in russ.)}  Uspekhi Mat.~Nauk 49(1994), no.~2, 159-160.
\smallskip

\item{[11]} Manuilov,~V.~M.:
{\it Diagonalization of compact operators on Hilbert modules over W*-algebras
of finite type (in engl.)}  submitted to Annals Global Anal. Geom., 1994.

\item{[12]} Murphy,~G.~J.:
{\it Diagonality in $C^*$-algebras.} Math.~Zeitschr. 199(1988), 199-229.
\smallskip

\item{[13]} Paschke,~W.~L.:
{\it Inner product modules over B*-algebras.} Trans.~Amer.~Math.~Soc.
182(1973), 443- 468.
\smallskip

\item{[14]} Wegge-Olsen,~N.~E.:
{\it K-theory and C*-algebras - a friendly approach.} Oxford University Press,
Oxford-New York-Tokyo, 1993.
\smallskip

\item{[15]} Zhang,~S.:
{\it Diagonalizing projections in the multiplier algebras and matrices
over a C*-algebra.} Pacific J.~Math. 145(1990), 181-200.
\smallskip

\item{[16]} Zhang,~S.:
{\it $K_1$-groups, quasidiagonality and interpolation by multiplier
projections.} Trans. Amer.~Math.~Soc. 325(1991), 793-818.
\vskip15mm

Received on June 28, 1994; revised on November 4, 1994.

To appear in {\it Zeitschrift f\"ur Analysis und ihre Anwendungen (ZAA)}
v. {\bf 14}(1995), no. 1, 33-41.

\bye